\documentclass[12pt]{article}

\makeatletter
\renewcommand{\@oddfoot}{\hfill \thepage}
\makeatother

\usepackage[english]{babel}
\usepackage{mathrsfs}
\usepackage{amssymb}
\usepackage{amsmath,amsthm,amsfonts}
\usepackage{cite}
\usepackage[dvipdf]{graphicx}
\usepackage{fancyhdr}
\usepackage[titletoc]{appendix} 

\begin{document}

\begin{center}
\Large{\bf Exact formula for the 2-marginal second moment function of the multidimensional 
symmetric\\ Markov random flight} 
\end{center}

\begin{center}
Alexander D. KOLESNIK\\
Institute of Mathematics and Computer Science\\
Moldova State University\\
Academy Street 5, Kishinev (Chi\c sin\u au) 2028, Moldova\\ 
E-Mail: alexander.kolesnik@math.md
\end{center}

\begin{abstract}
We consider the symmetric Markov random flight $\bold X(t), \; t>0,$ in the Euclidean space 
$\Bbb R^m, \; m\ge 3$, performed by a particle that moves in $\Bbb R^m$ with constant finite speed 
and changes its directions at Poisson-distributed random time instants by choosing the 
initial and each new direction at random according to the uniform distribution on the unit 
$(m-1)$-dimensional sphere. 
The 2-marginal second moment function $\mu_{(2,2,0,\dots,0)}(t), \; t>0,$ of $\bold X(t)$,  
corresponding to the multi-index $(2,2,0,\dots,0)$, is examined. An explicit formula for function 
$\mu_{(2,2,0,\dots,0)}(t)$ is obtained. This formula is also valid for all other 2-marginal second 
moment functions corresponding to any multi-indices of the form $(0,\dots,0,2,0,\dots,0,2,0,\dots,0)$. 
It is also shown that this moment function, under the standard Kac scaling condition, turns into 
the product of the variances of two coordinates of the $m$-dimensional homogeneous Brownian motion.  
\end{abstract}

\vskip 0.1cm


\bigskip 


\section{Introduction} 

Finite-velocity stochastic motions in the Euclidean spaces of various dimensions, also called {\it random flights}, have become an important field of modern stochastic analysis and statistical physics in recent decades. Their importance is determined by the fact that such stochastic motions generate transport processes that are the finite-velocity counterparts of the classical Einstein-Smoluchowski transport model with an infinite speed of propagation. While in Einstein-Smoluchowski model the transport is carried out by Brownian particles whose speed and intensity of changing the directions per unit of time are 
{\it infinite}, the random flights determine a transport by the particles moving with {\it finite} speed and changing their direction at random a {\it finite} number of times per unit of time. The finiteness of these basic parameters (speed and intensity of switching) implies essential difficulties in analyzing such stochastic motions in the Euclidean spaces of different dimensions. 

The Einstein-Smoluchowski transport model is governed by the parabolic heat equation with Laplace operator of respective dimension and, from this point of view, one can say that such processes with an infinite speed of propagation have a similar behaviour in the Euclidean spaces of arbitrary dimension. In contrast, the finite-velocity stochastic motions are driven by much more complicated equations and their behaviour is substantially different in the spaces of different dimensions.  

The stochastic motion at finite speed on the real line $\Bbb R^1$, whose evolution is controlled by a homogeneous Poisson process, was the first finite-velocity transport model appeared in literature. Apparently, V.A. Fock was the first to consider such stochastic motion in the framework of the  problem of describing the diffusion of a light ray in homogeneous environment. He has derived a 
one-dimensional hyperbolic telegraph equation, which this diffusion process was subject to 
(see \cite[Section 13, formula (130)]{fock}). But a true breakthrough occurred when the works by 
S. Goldstein \cite{gold} and M. Kac \cite{kac} appeared. Studying the process of stochastic motion at finite constant speed on the real line $\Bbb R^1$ with two alternating directions taking at Poisson-distributed random time instants, they have shown that the transition density of the motion is the fundamental solution 
(the Green's function) to the one-dimensional hyperbolic telegraph equation of second order with constant coefficients.  

Telegraph-type processes, also called {\it persistent random walks}, and their numerous generalizations have become a well-developed area of modern stochastic analysis, statistical physics, transport processes with diverse applications in physics \cite{bras,giona1,giona2,giona3,giona4}, relativity theory 
\cite{cane1,cane2, giona}, cosmology and astrophysics \cite{broad1,broad2,reim}, radioactive transport 
\cite{eon}, quantum physics \cite{gzyl}, financial modeling \cite{kolrat} and some other fields of science, technology and engineering. At present, this area of research numbers more than hundred of works, including several monographs. A vast bibliography on this subject can be found in \cite{mas1,mas2} and in the recent monographs \cite{kol1,kolrat}. In the last decade, mostly in physical literature, some properties of such finite-velocity stochastic motions have been studied in the framework of the run-and-tumble theory 
\cite{angel,dhar,mall,zhang}.  

Among the great variety of works devoted to the finite-velocity stochastic motions, the overwhelming  majority deal with one-dimensional models. As to the multidimensional counterparts of the telegraph-type  
processes is concerned, the palette of obtained results looks much poorer than in the one-dimensional case.  As noted above, this fact is explained by the great difficulties in analyzing such multidimensional finite-velocity stochastic motions. The natural multidimensional generalization of the Goldstein-Kac telegraph process is represented by the stochastic motion of a particle moving at finite constant speed in the Euclidean space $\Bbb R^m, \; m\ge 2$, and changing its directions at Poisson-distributed random time instants by choosing them at random according to the uniform probability law on the unit $(m-1)$-dimensional sphere. Such stochastic motion is referred to as the {\it Markov random flight}. 

When studying any stochastic process, its distribution is undoubtedly the most desirable goal. As to the 
Markov random flights is concerned, their distributions were obtained only in a few Euclidean spaces of low even dimensions. The explicit transition density of the symmetric Markov random flight in the Euclidean 
plane $\Bbb R^2$ was obtained by different methods in \cite{kol5,mas,sta1,sta2} (see also 
\cite[formula (5.2.5) ]{kol1}).  The transition density of the symmetric Markov random flight in the Euclidean space $\Bbb R^4$ was presented in \cite{kol4} (see also \cite[formula (7.2.4) ]{kol1}). 
Amazingly, both these transition densities in the spaces  $\Bbb R^2$ and  $\Bbb R^4$ are expressed in terms of elementary (exponential) function. Finally, the transition density of the symmetric Markov random flight in the space $\Bbb R^6$ was obtained in \cite{kol2} (see also \cite[formula (8.2.6)]{kol1}) and it has the form of a fairly complicated series with respect to Gauss hypergeometric functions. In all other Euclidean spaces, the distributions of the Markov random flights, including that in the extremely important three-dimensional space $\Bbb R^3$, was not obtained so far. 

In this situation, when the distribution of the Markov random flight in arbitrary dimension is unknown, 
the urgent necessity appears to study its other characteristics and, first of all, the numerical ones. 
It is well known that moments are one of the most significant numerical characteristics of a one-dimensional random variable. Among them, the first and the second moments are of a special importance. The first moment, that is the expectation, is  the mean value of all weighted values of random variable. But much more important characteristic of random variable is its second moment, that is, the variance, representing the measure of dispersion of the values of random variable around its mean value. Many fundamental results in probability theory were proved under the assumption of existence of finite-valued variance. Moreover, many probability equalities and inequalities include a variance parameter. All this is also related to a one-dimensional stochastic process, whose moment function is one of the most important characteristics, which, for any fixed value of the time variable, produces the exact value of moment. 

The situation becomes much more complicated in the case of multidimensional stochastic processes. 
The moment function $\mu_{\bf q}(t), \; t>0,$ of a $m$-dimensional stochastic process is determined by 
its multi-index ${\bf q} = (q_1,q_2,\dots,q_m), \; m\ge 2$, where $q_i, \; i=1,2,\dots,m$, are arbitrary 
non-negative real numbers (for rigorous definition of the moment function of a $m$-dimensional stochastic process, see formulas (\ref{eq01}) and (\ref{eq2}) below). The multidimensional counterpart of the variance is the second moment function $\mu_{\bf 2}(t)=\mu_{(2,2,\dots,2)}(t), \; t>0,$ determined by the 
multi-index ${\bf 2}=(2,2,\dots,2)$ of length $m$. 

Moment functions of the symmetric Markov random flights in some Euclidean spaces have already been studied in literature. Moment function of the Goldstein-Kac telegraph process (which is the 
one-dimensional Markov random flight), was examined in \cite{iacus,koles} 
(see also \cite[Section 2.9]{kol1}). Using the known distributions of the Markov random flights in the Euclidean spaces $\Bbb R^2$ and $\Bbb R^4$, their moment functions were obtained in \cite{kol03} 
(see also \cite[Sections 5.7 and 7.5]{kol1}). But finding of the moment functions of the Markov 
random flights in other dimensions seems to be an impracticable problem. That is why it would be quite reasonable to concentrate the efforts on analyzing the most important moment functions, namely, on the 
first and second ones.  

Usually, finding of the first moment function is a more easy task, especially if a stochastic process has 
symmetric structure. However, obtaining of the second moment function of a $m$-dimensional stochastic process, in general, is a much more difficult analytical problem, especially if the distribution of the process is unknown. That is why, finding of the marginals of the second moment function is of a special interest. In this case, the term '$l$-marginal' means that the respective multi-index contains only $l$ 
non-zero elements ($l<m$). Surprisingly, in some cases of not too big integers $l$ and the symmetric structure of the stochastic process, one manages to obtain the exact formulas for such $l$-marginal 
second moment functions.  

Everything said above applies entirely to the multidimensional symmetric Markov random flights. In particular, it was proved in \cite[Theorem 3]{kol0} that the first moment function of the three-dimensional 
symmetric Markov random flight, corresponding to the multi-index ${\bf 1}=(1,1,1)$, is identically equal to zero (that should be expected due to the symmetry of the process), while for the second moment function determined by the multi-index ${\bf 2}=(2,2,2)$ an asymptotic formula for small values of time variable,  was given. Moreover, it was also shown (see \cite[formula (70)]{kol0}) that the $1$-marginal second moment function of the symmetric Markov random flight of arbitrary dimension $m\ge 3$, corresponding to the multi-index $(2,0,\dots,0)$ of length $m$, (as well as to any other multi-index of the form 
$(0,\dots,0,2,0,\dots,0)$), is given by the formula:  
\begin{equation}\label{eq2900} 
\mu_{(2,0,\dots,0)}(t) = \frac{2}{m} \; \frac{c^2}{\lambda^2} \left( e^{-\lambda t} + \lambda t - 1 \right), 
\qquad m\ge 3, \quad t>0,
\end{equation} 
where $c$ and $\lambda$ are the speed of motion and intensity of switching the directions, respectively. 

In particular, the $1$-marginal second moment function corresponding to the multi-index $(2,0,0)$ of the 
three-dimensional symmetric Markov random flight, takes the form (see \cite[formula (69)]{kol0}):  
\begin{equation}\label{eq2910} 
\mu_{(2,0,0)}(t) = \frac{2}{3} \; \frac{c^2}{\lambda^2} \left( e^{-\lambda t} + \lambda t - 1 \right), 
\qquad \quad t>0.
\end{equation} 

Note that $1$-marginal second moment function (\ref{eq2910}) is resembled to that presented in 
\cite[formula (49)]{zhang} (for diffusive coefficient $D=0$ therein), however that formula contains wrong factor $2$, instead of correct factor $2/3$ in our formula (\ref{eq2910}). Moreover, the dimension 
parameter is absent in \cite[formula (49)]{zhang} at all, and this is impossible, because the characteristics of any multidimensional random variable or stochastic process (even with independent coordinates) must contain its dimension in that or another form. 

In this article we make the next important step in studying the second moment function of the multidimensional symmetric Markov random flight and obtain an explicit formula for its $2$-marginal 
second moment function valid in arbitrary dimension $m\ge 3$. The paper is organized as follows. 
In Section 2, for the reader's convenience, we recall some previous results that our further analysis 
is based on. In Sections 3 and 4, the explicit formulas are obtained for the coefficients in the polynomials, which the characteristic function of the Markov random flight is decomposed in, at the 
second and fourth powers of the norm of  inversion parameter. Using these formulas, as well as the 
well-known connection between the characteristic and moment functions, in Section 5, we prove a theorem yielding an exact relation for the $2$-marginal second moment function of the symmetric Markov random 
flight in arbitrary dimension $m\ge 3$. This exact formula is the main result of the article. Finally, 
in Section 6, we give some conclusions and final remarks concerning the further generalizations 
of the obtained results. In Appendix, two auxiliary lemmas are proved, which are used in the derivation 
of the main result.

\section{Preliminaries} 

The $m$-dimensional symmetric Markov random flight $\bold X(t)=(X_1(t),\dots,X_m(t)), \; t>0,$ is represented by the stochastic motion of a particle moving with constant finite speed $c>0$ in the Euclidean space $\Bbb R^m, \; m\ge 3,$ and taking on its initial and each new directions at $\lambda$-Poisson ($\lambda>0$) distributed random time instants by choosing them at random according to the uniform distribution on the unit $(m-1)$-dimensional sphere $S_1^m = \left\{ \bold x=(x_1, \dots ,x_m)\in \Bbb R^m: \;
\Vert\bold x\Vert^2 = x_1^2+ \dots +x_m^2=1 \right\}$. 

At arbitrary time instant $t>0$, the distribution 
$\text{Pr}\{ \bold X(t)\in d\bold x \}$ is concentrated in the closed ball 
$$\bold B_{ct}^m = \left\{ \bold x=(x_1, \dots ,x_m)\in \Bbb R^m : \;
\Vert\bold x\Vert^2 = x_1^2+ \dots +x_m^2\le c^2t^2 \right\} .$$
The singular component of the distribution is concentrated on the boundary 
$$S_{ct}^m =\partial\bold B_{ct}^m = \left\{ \bold x=(x_1,
\dots ,x_m)\in \Bbb R^m: \; \Vert\bold x\Vert^2 = x_1^2+ \dots +x_m^2=c^2t^2 \right\}$$
of the ball $\bold B_{ct}^m$, while its absolutely continuous part is concentrated in the interior 
$$\text{int} \; \bold B_{ct}^m = \left\{ \bold x=(x_1,
\dots ,x_m)\in \Bbb R^m: \; \Vert\bold x\Vert^2 = x_1^2+ \dots +x_m^2<c^2t^2 \right\}$$
of this ball. 

Consider the characteristic function of the Markov random flight $\bold X(t)$ defined by 
\begin{equation}\label{eq1}
H(\boldsymbol\alpha,t) = \Bbb E \left\{ e^{i\langle\boldsymbol\alpha,\bold X(t)\rangle} \right\} ,  
\qquad t>0,
\end{equation}
where $\Bbb E$ means the expectation, $\boldsymbol\alpha=(\alpha_1, \dots ,\alpha_m) \in\Bbb R^m$ is 
the real-valued $m$-dimensional vector of inversion parameters and 
$\langle\boldsymbol\alpha,\bold X(t)\rangle$ 
is the inner product of the vectors $\boldsymbol\alpha$ and $\bold X(t)$.

The connection between the characteristic function of a stochastic process and its moment function 
is well known. We will use this connection in order to obtain some important relations for the 
marginals of the moment function of the $m$-dimensional symmetric Markov random flight 
$\bold X(t), \; t>0$. 

Let $\bold q=(q_1,\dots, q_m)$ be a multi-index. Since at arbitrary fixed time $t>0$ the Markov random flight $\bold X(t) = (X_1(t),\dots,X_m(t))$ is concentrated in the ball $\bold B_{ct}^m$, this process, 
as well as its coordinates $X_j(t), \; j=1,\dots,m$, are bounded and, therefore, for arbitrary positive integer $k$, the condition $\Bbb E |X_j(t)|^k<\infty$ fulfills for all $j=1,\dots,m$. Then, for arbitrary fixed $t>0$, there exist the mixed moments   
\begin{equation}\label{eq01}
\mu_{\bold q}(t) = \mu_{(q_1,\dots,q_m)}(t) = \Bbb E X_1^{q_1}(t)\dots X_m^{q_m}(t) , 
\qquad q_j\ge 0, \quad q_1+\dots +q_m \le k,
\end{equation} 
of the $m$-dimensional symmetric Markov random flight $\bold X(t)=(X_1(t),\dots,X_m(t))$ given by the formula: 
\begin{equation}\label{eq2} 
\mu_{\bold q}(t) = \mu_{(q_1,\dots,q_m)}(t) = (-i)^{q_1+\dots +q_m} \; 
\frac{\partial^{q_1+\dots +q_m}}{\partial\alpha_1^{q_1} \dots \partial\alpha_m^{q_m}} \; 
H(\alpha_1,\dots,\alpha_m; t) \biggr|_{\alpha_1= \dots =\alpha_m=0} . 
\end{equation}
This general formula (\ref{eq2}) will be used for obtaining the 2-marginals of the second moment function 
of the Markov random flight $\bold X(t)$. 

Two series representations for the characteristic function (\ref{eq1}) of the $m$-dimensional symmetric Markov random flight $\bold X(t)$ were obtained in recent article \cite{kol0}. The first representation is given in the form of a series with respect to Bessel functions with variable indices 
(see \cite[Theorem 1]{kol0}). However, for our purposes it is more convenient to use the second representation of the characteristic function in the form of a series with respect to the powers of time variable. This representation is given by the relation (see \cite[Theorem 2]{kol0}):
\begin{equation}\label{eq3} 
H(\boldsymbol\alpha,t) = e^{-\lambda t} \sum_{n=0}^{\infty} 
\frac{\gamma_{n+1}(\Vert\boldsymbol\alpha\Vert)}{n!} \; t^n ,
\end{equation} 
where the coefficients $\gamma_n=\gamma_n(\Vert\boldsymbol\alpha\Vert)$ are given by the recurrent relation: 
\begin{equation}\label{eq4} 
\gamma_1 = 1,  \qquad \gamma_n = \theta_n + \lambda \sum_{k=1}^{n-1} 
\gamma_{n-k} \; \theta_k , \qquad n\ge 2, 
\end{equation} 
and  
\begin{equation}\label{eq5}
\theta_n = \left\{ 
\aligned 0, \qquad & \quad \text{if} \;\; n=2k, \\
         (-1)^k \; \frac{\left( \frac{1}{2} \right)_k}{\left( \frac{m}{2} \right)_k} \; 
				(c\Vert\boldsymbol\alpha\Vert)^{2k} , 
				& \quad \text{if} \;\; n=2k+1 , 
				\endaligned \right. \qquad k=0,1,2,\dots , 
\end{equation}
where 
$$(a)_n=a(a+1)\dots (a+n-1)=\frac{\Gamma (a+n)}{\Gamma (a)}$$ 
is the Pochhammer symbol.

Since all the coefficients (\ref{eq5}) with even indices are zeros, then recurrent relation (\ref{eq4}) 
can be split for even and odd indices as follows (\cite[Remark 6]{kol0}):
\begin{equation}\label{eq6} 
\gamma_{2n} = \lambda \sum_{k=0}^{n-1} \gamma_{2n-2k-1} \; \theta_{2k+1} , \qquad n\ge 1, 
\end{equation}
\begin{equation}\label{eq7} 
\gamma_1=1, \qquad \gamma_{2n+1} = \theta_{2n+1} +\lambda \sum_{k=0}^{n-1} \gamma_{2n-2k} \; \theta_{2k+1} , \qquad n\ge 1, 
\end{equation}
where, according to (\ref{eq5}), 
\begin{equation}\label{eq8}
\theta_1 = 1, \qquad  \theta_{2k+1} = (-1)^k \; \frac{\left( \frac{1}{2} \right)_k}{\left( \frac{m}{2} 
\right)_k} \; (c\Vert\boldsymbol\alpha\Vert)^{2k} , \qquad k\ge 1. 
\end{equation}

The eight coefficients $\gamma_n$ are given by the formulas (see \cite[relations (49)]{kol0}): 
\begin{equation}\label{eq9} 
\aligned 
& \gamma_0 = 0, \qquad \gamma_1 = 1, \qquad \gamma_2 = \lambda, \qquad 
\gamma_3 = \lambda^2 - \frac{1}{m} \; (c\Vert\boldsymbol\alpha\Vert)^2 , \\  
& \gamma_4 = \lambda^3 - \lambda \; \frac{2}{m} \; (c\Vert\boldsymbol\alpha\Vert)^2 , \qquad 
\gamma_5 = \lambda^4 - \lambda^2 \; \frac{3}{m} \; (c\Vert\boldsymbol\alpha\Vert)^2 
+ \frac{3}{m(m+2)} \; (c\Vert\boldsymbol\alpha\Vert)^4 \\
& \gamma_6 = \lambda^5 - \lambda^3 \; \frac{4}{m} \; (c\Vert\boldsymbol\alpha\Vert)^2 
+ \lambda \; \frac{7m+2}{m^2 (m+2)} \; (c\Vert\boldsymbol\alpha\Vert)^4 \\ 
& \gamma_7 = \lambda^6 - \lambda^4 \; \frac{5}{m} \; (c\Vert\boldsymbol\alpha\Vert)^2 
+ \lambda^2 \; \frac{12m+6}{m^2 (m+2)} \; (c\Vert\boldsymbol\alpha\Vert)^4  
- \frac{15}{m(m+2)(m+4)} \; (c\Vert\boldsymbol\alpha\Vert)^6 \\ 
& \gamma_8 = \lambda^7 - \lambda^5 \; \frac{6}{m} \; (c\Vert\boldsymbol\alpha\Vert)^2  
+ \lambda^3 \; \frac{18m+12}{m^2 (m+2)} \; (c\Vert\boldsymbol\alpha\Vert)^4  
- \lambda \; \frac{36m+24}{m^2 (m+2)(m+4)} \; (c\Vert\boldsymbol\alpha\Vert)^6 \\ 
\endaligned 
\end{equation}

One can see that, for arbitrary $n\ge 1$, the coefficient 
$\gamma_n = \gamma_n(\Vert\boldsymbol\alpha\Vert)$ is a polynomial of the variable 
$\Vert\boldsymbol\alpha\Vert$ of the power $2\left[ \frac{n-1}{2} \right], \; n\ge 1,$ or, equivalently, 
a polynomial of the variable $(\alpha_1^2+ \dots +\alpha_m^2)$ of the power 
$\left[ \frac{n-1}{2} \right], \; n\ge 1$, where $[ \; \cdot \; ]$ means the integer part of a number. 
Therefore, each polynomial $\gamma_n = \gamma_n(\Vert\boldsymbol\alpha\Vert)$ has the form: 
\begin{equation}\label{PolynomGamma} 
\gamma_n = \gamma_n(\Vert\boldsymbol\alpha\Vert) = \sum_{k=0}^{[(n-1)/2]} 
K_{\Vert\boldsymbol\alpha\Vert^{2k}}(\gamma_n) \; \Vert\boldsymbol\alpha\Vert^{2k} , 
\qquad n\ge 1 , 
\end{equation}
where $K_{\Vert\boldsymbol\alpha\Vert^{2k}}(\gamma_n)$ are the coefficients at the term 
$\Vert\boldsymbol\alpha\Vert^{2k}$ in the polynomials $\gamma_n = \gamma_n(\Vert\boldsymbol\alpha\Vert)$. 
Obviously, $K_{\Vert\boldsymbol\alpha\Vert^{2\cdot 0}}(\gamma_n) = \lambda^{n-1}$ (see (\ref{eq9})), and, 
for arbitrary $n\ge 2$, the coefficients $K_{\Vert\boldsymbol\alpha\Vert^{2k}}(\gamma_n)$ depend on the parameters $c$ and $\lambda$, as well as on the dimension $m$. 

The polynomials $\gamma_n = \gamma_n(\Vert\boldsymbol\alpha\Vert), \; n\ge 1,$ play a key role in obtaining the closed-form expressions for the characteristic and moment functions of the 
$m$-dimensional symmetric Markov random flight $\bold X(t)=(X_1(t),\dots,X_m(t)), \; t>0$. The crucial 
point is to derive a general formula for the coefficients $K_{\Vert\boldsymbol\alpha\Vert^{2k}}(\gamma_n)$ at the term $\Vert\boldsymbol\alpha\Vert^{2k}$ in arbitrary polynomials 
$\gamma_n = \gamma_n(\Vert\boldsymbol\alpha\Vert), \; n\ge 1$. This, however, is a very difficult (and, maybe, impracticable) analytical problem. Nevertheless, in the next sections, we will be able to obtain exact formulas for the coefficients $K_{\Vert\boldsymbol\alpha\Vert^2}(\gamma_n)$ and 
$K_{\Vert\boldsymbol\alpha\Vert^4}(\gamma_n)$ (corresponding to the numbers $k=1$ and $k=2$ in series 
(\ref{PolynomGamma})) at the terms $\Vert\boldsymbol\alpha\Vert^2$ and 
$\Vert\boldsymbol\alpha\Vert^4$, respectively, in arbitrary polynomials 
$\gamma_n = \gamma_n(\Vert\boldsymbol\alpha\Vert), \; n\ge 1$.
This enables us to derive an explicit relation for the 2-marginal second moment function and this result 
is an important step in studying the basic characteristics of the multidimensional Markov random flight 
$\bold X(t)$.

\section{Explicit formula for the coefficients at the term $\Vert\boldsymbol\alpha\Vert^2$ in the polynomials $\gamma_n(\Vert\boldsymbol\alpha\Vert)$}

In this section we will derive an exact formula for the coefficients 
$K_{\Vert\boldsymbol\alpha\Vert^2}(\gamma_n)$ in the polynomials 
$\gamma_n=\gamma_n(\Vert\boldsymbol\alpha\Vert)$ related to the term $\Vert\boldsymbol\alpha\Vert^2$, 
which will be used in the proof of the main result. This formula has already been noted in 
\cite[page 20]{kol0} with some hints of its deriving, but no proof was given therein. 

In the following proposition we fill this gap and give a rigorous derivation of the formula related to 
the coefficients mentioned above. 

\bigskip 

{\bf Proposition 1.} {\it For arbitrary dimension $m\ge 3$ and for arbitrary $n\ge 3$, the coefficients 
at the term $\Vert\boldsymbol\alpha\Vert^2$ in the polynomials 
$\gamma_n=\gamma_n(\Vert\boldsymbol\alpha\Vert)$ of the $m$-dimensional symmetric Markov random flight 
$\bold X(t)$ are given by the formula:} 
\begin{equation}\label{CoeffA2} 
K_{\Vert\boldsymbol\alpha\Vert^2}(\gamma_n) = - \frac{n-2}{m} \; c^2 \; \lambda^{n-3} , 
\qquad n\ge 3, \quad m\ge 3 .
\end{equation}

\begin{proof}

We will prove (\ref{CoeffA2}) by induction. For $n=3$ and $n=4$, formula (\ref{CoeffA2}) yields: 
$$K_{\Vert\boldsymbol\alpha\Vert^2}(\gamma_3) = - \frac{1}{m} \; c^2 , \qquad 
K_{\Vert\boldsymbol\alpha\Vert^2}(\gamma_4) = - \frac{2}{m} \; c^2 \; \lambda ,$$ 
and this exactly coincides with respective coefficients in (\ref{eq9}).

Since recurrent relations (\ref{eq6}) and (\ref{eq7}) are somewhat different, we will prove (\ref{eq9}) separately for even numbers $n\ge 6$ and odd numbers $n\ge 5$.

$\bullet$ {\bf The case of even $n\ge 6$}. Suppose that formula (\ref{CoeffA2}) is true for all the numbers 
$3,4,\dots,(2n-1)$. Our aim is to prove it for the number $2n$. Since $n$ is even, that is, 
$n=2k, \; k\ge 2$, and since the coefficient $\theta_{2k+1}$ contains the factor  
$\Vert\boldsymbol\alpha\Vert^{2k}$ (see (\ref{eq8})), then we should take into account only those $k$, 
which satisfy the inequality $2k\le 2$, that is only $k=0$ and $k=1$. 

In this case, according to recurrent relation (\ref{eq6}), we have: 
\begin{equation}\label{eq10} 
K_{\Vert\boldsymbol\alpha\Vert^2}(\gamma_{2n}) = 
\lambda \left[ K_{\Vert\boldsymbol\alpha\Vert^2}(\gamma_{2n-1} \; \theta_1) + 
K_{\Vert\boldsymbol\alpha\Vert^2}(\gamma_{2n-3} \; \theta_3) \right] . 
\end{equation}
Since $\theta_1=1$ and according to induction assumption, we get: 
\begin{equation}\label{eq11} 
K_{\Vert\boldsymbol\alpha\Vert^2}(\gamma_{2n-1} \; \theta_1) = 
K_{\Vert\boldsymbol\alpha\Vert^2}(\gamma_{2n-1}) = - \frac{2n-3}{m} \; c^2 \; \lambda^{2n-4} .
\end{equation} 

Similarly, in view of (\ref{eq8}) and according to induction assumption, we have: 
\begin{equation}\label{eq12} 
K_{\Vert\boldsymbol\alpha\Vert^2}(\gamma_{2n-3} \; \theta_3) = 
K_{\Vert\boldsymbol\alpha\Vert^0}(\gamma_{2n-3}) \;  
K_{\Vert\boldsymbol\alpha\Vert^2}(\theta_3) =
- \lambda^{2n-4} \; \frac{\left(\frac{1}{2}\right)_1}{\left(\frac{m}{2}\right)_1} \; c^2 
= - \frac{1}{m} \; \lambda^{2n-4} \; c^2 ,  
\end{equation} 
where we have used the fact that 
\begin{equation}\label{Poch1} 
\frac{\left(\frac{1}{2}\right)_1}{\left(\frac{m}{2}\right)_1} = \frac{1}{m} .
\end{equation} 

Substituting (\ref{eq11}) and (\ref{eq12}) into (\ref{eq10}), after some simple calculations, 
we finally obtain: 
$$K_{\Vert\boldsymbol\alpha\Vert^2}(\gamma_{2n}) = - \frac{(2n)-2}{m} \; c^2 \; \lambda^{(2n)-3} ,$$
proving (\ref{CoeffA2}) for even $n\ge 4$.

\bigskip 

$\bullet$ {\bf The case of odd $n\ge 5$}. Suppose that formula (\ref{CoeffA2}) is true for all the numbers 
$3,4,\dots,(2n)$. Our aim is to prove it for the number $(2n+1)$. For the same reason as above, we should take into account only the two values of variable index $k$, namely, $k=0$ and $k=1$. 
Then, according to (\ref{eq7}), we have: 
\begin{equation}\label{eq13} 
K_{\Vert\boldsymbol\alpha\Vert^2}(\gamma_{2n+1}) = K_{\Vert\boldsymbol\alpha\Vert^2}(\theta_{2n+1}) + 
\lambda \left[ K_{\Vert\boldsymbol\alpha\Vert^2}(\gamma_{2n} \; \theta_1) + 
K_{\Vert\boldsymbol\alpha\Vert^2}(\gamma_{2n-2} \; \theta_3) \right] . 
\end{equation}
It is obvious that for any $n\ge 3$, 
\begin{equation}\label{eq14} 
K_{\Vert\boldsymbol\alpha\Vert^2}(\theta_{2n+1}) = 0 .
\end{equation}

Then, taking into account that $\theta_1=1$ and by induction assumption, we get: 
\begin{equation}\label{eq15}
K_{\Vert\boldsymbol\alpha\Vert^2}(\gamma_{2n} \; \theta_1) 
= K_{\Vert\boldsymbol\alpha\Vert^2}(\gamma_{2n}) = - \frac{2n-2}{m} \; c^2 \; \lambda^{2n-3} .
\end{equation}
Similarly, in view of (\ref{eq8}) and according to induction assumption, we have: 
\begin{equation}\label{eq16} 
K_{\Vert\boldsymbol\alpha\Vert^2}(\gamma_{2n-2} \; \theta_3) = 
K_{\Vert\boldsymbol\alpha\Vert^0}(\gamma_{2n-2}) \;  
K_{\Vert\boldsymbol\alpha\Vert^2}(\theta_3) =
- \lambda^{2n-3} \; \frac{\left(\frac{1}{2}\right)_1}{\left(\frac{m}{2}\right)_1} \; c^2 
= - \frac{1}{m} \; c^2 \; \lambda^{2n-3} . 
\end{equation}

Substituting (\ref{eq15}) and (\ref{eq16}) into (\ref{eq13}) and taking into account (\ref{eq14}), 
after some simple calculations, we finally obtain: 
$$K_{\Vert\boldsymbol\alpha\Vert^2}(\gamma_{2n+1}) = - \frac{(2n+1)-2}{m} \; c^2 \; \lambda^{(2n+1)-3} ,$$
proving (\ref{CoeffA2}) for odd $n\ge 5$.

The proposition is thus completely proved. 
\end{proof}

\section{Explicit formula for the coefficients at the term $\Vert\boldsymbol\alpha\Vert^4$ in the polynomials $\gamma_n(\Vert\boldsymbol\alpha\Vert)$} 

In this section we will obtain an explicit formula for the coefficients 
$K_{\Vert\boldsymbol\alpha\Vert^4}(\gamma_n)$ in the polynomials 
$\gamma_n=\gamma_n(\Vert\boldsymbol\alpha\Vert)$ related to the term $\Vert\boldsymbol\alpha\Vert^4$, 
which will be used in the proof of the main result. 

\bigskip

{\bf Proposition 2.} {\it For arbitrary dimension $m\ge 3$ and for arbitrary $n\ge 5$, the coefficients 
at the term $\Vert\boldsymbol\alpha\Vert^4$ in the polynomials 
$\gamma_n=\gamma_n(\Vert\boldsymbol\alpha\Vert)$ of the $m$-dimensional symmetric Markov random flight 
$\bold X(t)$ are given by the formula:} 
\begin{equation}\label{CoeffA4} 
K_{\Vert\boldsymbol\alpha\Vert^4}(\gamma_n) = \frac{c^4 \; \lambda^{n-5} \; (n-4)}{m^2 \; (m+2)} \; 
\left[ \frac{n+1}{2} \; m + (n-5) \right] ,  \qquad n\ge 5, \quad m\ge 3 .
\end{equation}

\begin{proof}

We will prove (\ref{CoeffA4}) by induction. For $n=5$ and $n=6$, formula (\ref{CoeffA4}) yields: 
$$K_{\Vert\boldsymbol\alpha\Vert^4}(\gamma_5) = \frac{3}{m \; (m+2)} \; c^4 , 
\qquad 
K_{\Vert\boldsymbol\alpha\Vert^4}(\gamma_6) = \frac{7m+2}{m^2 \; (m+2)} \; c^4 \; \lambda ,$$ 
and this exactly coincides with respective coefficients in (\ref{eq9}).

As in the proof of Proposition 1, we consider separately the cases of even numbers $n\ge 8$ and 
odd numbers $n\ge 7$.

$\bullet$ {\bf The case of even $n\ge 8$}. Suppose that formula (\ref{CoeffA4}) is true for all the numbers 
$7,8,\dots,(2n-1)$. Our aim is to prove it for the number $2n$. Since $n$ is even, that is, 
$n=2k, \; k\ge 3$, and since the coefficient $\theta_{2k+1}$ contains the factor  
$\Vert\boldsymbol\alpha\Vert^{2k}$ (see (\ref{eq8})), then we should take into account only those $k$, 
which satisfy the inequality $2k\le 4$, that is, only $k=0, \; k=1$ and $k=2$. 

Then, according to recurrent relation (\ref{eq6}), we have: 
\begin{equation}\label{eq17} 
K_{\Vert\boldsymbol\alpha\Vert^4}(\gamma_{2n}) = 
\lambda \left[ K_{\Vert\boldsymbol\alpha\Vert^4}(\gamma_{2n-1} \; \theta_1) + 
K_{\Vert\boldsymbol\alpha\Vert^4}(\gamma_{2n-3} \; \theta_3) 
+ K_{\Vert\boldsymbol\alpha\Vert^4}(\gamma_{2n-5} \; \theta_5) \right] . 
\end{equation} 

Let us evaluate separately the terms in (\ref{eq17}). Since $\theta_1=1$ and by induction assumption, 
for the first term in square brackets of (\ref{eq17}), we get : 
\begin{equation}\label{eq18} 
K_{\Vert\boldsymbol\alpha\Vert^4}(\gamma_{2n-1} \; \theta_1) = 
K_{\Vert\boldsymbol\alpha\Vert^4}(\gamma_{2n-1}) = \frac{c^4 \; \lambda^{2n-6} \; (2n-5)}{m^2 \; (m+2)} 
\; \bigl[ n m + (2n-6) \bigr] . 
\end{equation} 

According to (\ref{eq8}) and formula (\ref{CoeffA2}) of Proposition 1, for the second term in 
square brackets of (\ref{eq17}), we have:
\begin{equation}\label{eq19} 
\aligned 
K_{\Vert\boldsymbol\alpha\Vert^4}(\gamma_{2n-3} \; \theta_3) & = 
K_{\Vert\boldsymbol\alpha\Vert^2}(\gamma_{2n-3}) \;  
K_{\Vert\boldsymbol\alpha\Vert^2}(\theta_3) \\ 
& = \biggl[- \frac{2n-5}{m} \; c^2 \; \lambda^{2n-6} \biggr] \; 
\biggl[ - \frac{\left( \frac{1}{2} \right)_1}{\left( \frac{m}{2} \right)_1}  \; c^2 \biggr] \\ 
& = \frac{2n-5}{m^2} \; c^4 \; \lambda^{2n-6} . 
\endaligned
\end{equation} 
where we have again used (\ref{Poch1}).

Finally, for the third term in square brackets of (\ref{eq17}), we get:
\begin{equation}\label{eq20} 
K_{\Vert\boldsymbol\alpha\Vert^4}(\gamma_{2n-5} \; \theta_5) 
= K_{\Vert\boldsymbol\alpha\Vert^0}(\gamma_{2n-5}) \;  
K_{\Vert\boldsymbol\alpha\Vert^4}(\theta_3) 
= \lambda^{2n-6} \; \frac{\left( \frac{1}{2} \right)_2}{\left( \frac{m}{2} \right)_2}  \; c^4 = 
\frac{3}{m \; (m+2)} \; c^4 \; \lambda^{2n-6} , 
\end{equation} 
where we have used the fact that 
\begin{equation}\label{Poch2}
\frac{\left( \frac{1}{2} \right)_2}{\left( \frac{m}{2} \right)_2} = \frac{3}{m \; (m+2)} . 
\end{equation} 

Substituting (\ref{eq18}), (\ref{eq19}) and (\ref{eq20}) into (\ref{eq17}), after some simple calculations, 
we finally obtain: 
\begin{equation}\label{eq21} 
K_{\Vert\boldsymbol\alpha\Vert^4}(\gamma_{2n}) = 
\frac{c^4 \; \lambda^{(2n)-5} \; ((2n)-4)}{m^2 \; (m+2)} \; \left[ \frac{(2n)+1}{2} m + ((2n)-5) \right] , 
\end{equation} 
proving (\ref{CoeffA4}) for even $n\ge 6$.

\bigskip

$\bullet$ {\bf The case of odd $n\ge 7$}. Suppose that formula (\ref{CoeffA4}) is true for all the numbers 
$7,8,\dots,(2n)$. Our aim is to prove it for the number $(2n+1)$. For the same reason as above, we should consider only the numbers $k=0, \; k=1$ and $k=2$ in the sum of recurrent relation (\ref{eq7}). 
Then, according to (\ref{eq7}), we have:  
\begin{equation}\label{eq22} 
K_{\Vert\boldsymbol\alpha\Vert^4}(\gamma_{2n+1}) = K_{\Vert\boldsymbol\alpha\Vert^4}(\theta_{2n+1})
+ \lambda \left[ K_{\Vert\boldsymbol\alpha\Vert^4}(\gamma_{2n} \; \theta_1) + 
K_{\Vert\boldsymbol\alpha\Vert^4}(\gamma_{2n-2} \; \theta_3) 
+ K_{\Vert\boldsymbol\alpha\Vert^4}(\gamma_{2n-4} \; \theta_5) \right] . 
\end{equation} 

Obviously, for any $n\ge 3$, 
\begin{equation}\label{eq23} 
K_{\Vert\boldsymbol\alpha\Vert^4}(\theta_{2n+1}) =0 .
\end{equation}

Let us evaluate separately the terms in square brackets of (\ref{eq22}). Since $\theta_1=1$ and by 
induction assumption, for the first term, we get : 
\begin{equation}\label{eq24} 
K_{\Vert\boldsymbol\alpha\Vert^4}(\gamma_{2n} \; \theta_1) = 
K_{\Vert\boldsymbol\alpha\Vert^4}(\gamma_{2n}) = \frac{c^4 \; \lambda^{2n-5} \; (2n-4)}{m^2 \; (m+2)} 
\; \biggl[ \frac{2n+1}{2} \; m + (2n-5) \biggr] . 
\end{equation} 

In view of formula (\ref{CoeffA2}) of Proposition 1, for the second term in square brackets 
of (\ref{eq22}), we have: 
\begin{equation}\label{eq25} 
\aligned 
K_{\Vert\boldsymbol\alpha\Vert^4}(\gamma_{2n-2} \; \theta_3) & = 
K_{\Vert\boldsymbol\alpha\Vert^2}(\gamma_{2n-2}) \;  
K_{\Vert\boldsymbol\alpha\Vert^2}(\theta_3) \\ 
& = \left[ -\frac{2n-4}{m} \; c^2 \; \lambda^{2n-5} \right] 
\left[ - c^2 \; \frac{\left(\frac{1}{2}\right)_1}{\left(\frac{m}{2}\right)_1} \right] \\ 
& = \frac{2n-4}{m^2} \; c^4 \; \lambda^{2n-5} , 
\endaligned
\end{equation}
where we have again used (\ref{Poch1}).

For the third term in square brackets of (\ref{eq22}), we have: 
\begin{equation}\label{eq26} 
K_{\Vert\boldsymbol\alpha\Vert^4}(\gamma_{2n-4} \; \theta_5) = 
K_{\Vert\boldsymbol\alpha\Vert^0}(\gamma_{2n-4}) \;  
K_{\Vert\boldsymbol\alpha\Vert^4}(\theta_5) = 
\lambda^{2n-5} \; \frac{\left(\frac{1}{2}\right)_2}{\left(\frac{m}{2}\right)_2} \; c^4 
= \frac{3}{m \; (m+2)} \; c^4 \; \lambda^{2n-5} ,
\end{equation}
where we have used (\ref{Poch2}).

Substituting (\ref{eq24}), (\ref{eq25}) and (\ref{eq26}) into (\ref{eq22}), and taking into account 
(\ref{eq23}), after some calculations, we finally obtain: 
\begin{equation}\label{eq27} 
K_{\Vert\boldsymbol\alpha\Vert^4}(\gamma_{2n+1}) = 
\frac{c^4 \; \lambda^{(2n+1)-5} \; ((2n+1)-4)}{m^2 \; (m+2)} \; 
\left[ \frac{(2n+1)+1}{2} m + ((2n+1)-5) \right] , 
\end{equation} 
proving (\ref{CoeffA4}) for odd $n\ge 7$.
The proposition is thus completely proved. 
\end{proof}

\section{Exact formula for the moment function $\mu_{(2,2,0,\dots,0)}(t)$} 

In the previous sections, we have obtained explicit formulas for the coefficients at the terms 
$\Vert\boldsymbol\alpha\Vert^2$ and $\Vert\boldsymbol\alpha\Vert^4$ in the polynomials 
$\gamma_n = \gamma_n(\Vert\boldsymbol\alpha\Vert)$ for arbitrary $n\ge 3$. With this in hand, 
we can now prove our main result related to the 2-marginal second moment function 
$\mu_{(2,2,0,\dots,0)}(t)$ of the symmetric Markov random flight $\bold X(t)$ in the $m$-dimensional Euclidean space $\Bbb R^m$ of arbitrary dimension $m\ge 3$. Amazingly, this 2-marginal second moment 
function $\mu_{(2,2,0,\dots,0)}(t)$ (as well as the 2-marginal second  moment functions 
corresponding to any multi-indices of the form $(0,\dots,0,2,0,\dots,0,2, 0,\dots,0)$) can be obtained 
in an explicit form. This main result of the article is given by the following theorem. 

\bigskip 

{\bf Theorem 1.} {\it For arbitrary dimension $m\ge 3$, the 2-marginal second moment function 
$\mu_{(2,2,0,\dots,0)}(t), \; t>0,$ of the $m$-dimensional symmetric Markov random 
flight $\bold X(t)$ is given by the formula:} 
\begin{equation}\label{SecondMomentFunct} 
\aligned 
\mu_{(2,2,0,\dots,0)}(t) & = 
\frac{8 \; c^4}{m^2 \; (m+2) \; \lambda^4} \biggl\{   
m \left[ e^{-\lambda t} \bigl( \lambda^2t^2 +3\lambda t + 3 \bigr) + 
\frac{\lambda^2t^2}{2} - 3 \right] \\
& \hskip 4cm + (\lambda^2t^2 + 12) (1-e^{-\lambda t}) - 6\lambda t (1+e^{-\lambda t}) \biggr\} . 
\endaligned
\end{equation}

\begin{proof}
In this particular case, general formula (\ref{eq2}) becomes: 
\begin{equation}\label{eq0191} 
\mu_{(2,2,0,\dots,0)}(t) = (-i)^4 \; \frac{\partial^4}{\partial\alpha_1^2 \; \partial\alpha_2^2} \; 
H(\alpha_1,\alpha_2,\dots,\alpha_m;t) \biggr|_{\alpha_1=\alpha_2= \dots =\alpha_m=0} .
\end{equation}

For the sake of brevity, let us introduce the operator: 
\begin{equation}\label{eq311} 
\mathcal A \; \cdot = 
\frac{\partial^4 \; \cdot}{\partial\alpha_1^2 \; \partial\alpha_2^2} \; 
\biggr|_{\alpha_1=\alpha_2=\dots=\alpha_m=0}   
\end{equation} 
Then, applying operator (\ref{eq311}) to the characteristic function given by series representation 
(\ref{eq3}) and taking into account that (see (\ref{eq9})) 
$\mathcal A \; \gamma_1=\mathcal A \; \gamma_2 =\mathcal A \; \gamma_3 =\mathcal A \; \gamma_4=0$, we have:   
\begin{equation}\label{eq201} 
\aligned 
\mathcal A \; H(\alpha_1,\alpha_2,\dots,\alpha_m; t) & = e^{-\lambda t} \sum_{k=4}^{\infty} 
\left[ \mathcal A \; \gamma_{k+1}(\alpha_1,\alpha_2,\dots,\alpha_m) \right] \; \frac{t^k}{k!} \\
& = e^{-\lambda t} \sum_{k=5}^{\infty} \frac{t^{k-1}}{(k-1)!}
\left[ \mathcal A \; \gamma_k(\alpha_1,\alpha_2,\dots,\alpha_m) \right]  .
\endaligned
\end{equation} 

It is easy to see that 
\begin{equation}\label{eq291} 
\mathcal A \; \Vert\boldsymbol\alpha\Vert^{2k} = \frac{\partial^4}{\partial\alpha_1^2 \; \partial\alpha_2^2} \; (\alpha_1^2+\alpha_2^2+\dots+\alpha_m^2)^k \; \biggr|_{\alpha_1=\alpha_2=\dots=\alpha_m=0} 
= \left\{ \aligned 
8 , \qquad & \text{if} \;\; k=2 ,\\ 
0, \qquad & \text{if} \;\; k\ge 3 , \endaligned \right.    
\end{equation} 
$$\boldsymbol\alpha = (\alpha_1,\alpha_2,\dots,\alpha_m)\in\Bbb R^m,$$ 
and, therefore, we should take into account only those terms of the polynomials 
$\gamma_k(\boldsymbol\alpha)$, which contain the factor 
$\Vert\boldsymbol\alpha\Vert^4 = (\alpha_1^2+\alpha_2^2+\dots+\alpha_m^2)^2$. 

Hence, according to formula (\ref{CoeffA4}) of Proposition 2 and taking into account (\ref{eq291}), 
we have: 
$$\mathcal A \; \gamma_k(\alpha_1,\alpha_2,\dots,\alpha_m) = 
\frac{8 \; c^4 \; \lambda^{k-5} \; (k-4)}{m^2 \; (m+2)} \; \left[ \frac{k+1}{2} \; m + (k-5) \right] , 
\qquad k\ge 5 .$$

Substituting this into (\ref{eq201}), we obtain: 
$$\aligned 
\mathcal A \; & H(\alpha_1,\alpha_2,\dots,\alpha_m; t) \\ 
& = 8 \; e^{-\lambda t} \sum_{k=5}^{\infty} \frac{t^{k-1}}{(k-1)!} \; 
\frac{c^4 \; \lambda^{k-5} \; (k-4)}{m^2 \; (m+2)} \; \left[ \frac{k+1}{2} \; m + (k-5) \right] \\ 
& = \frac{8 \; e^{-\lambda t} \; c^4}{m^2 \; (m+2)} \sum_{k=4}^{\infty} \frac{t^k}{k!} \; 
\lambda^{k-4} \; (k-3) \; \; \left[ \frac{k+2}{2} \; m + (k-4) \right] \\ 
& = \frac{8 \; e^{-\lambda t} \; c^4}{m^2 \; (m+2) \; \lambda^4} \sum_{k=4}^{\infty} 
\frac{(\lambda t)^k}{k!} \; \left[ \frac{m}{2} \; (k-3)(k+2) + (k-3)(k-4) \right] \\ 
& = \frac{8 \; e^{-\lambda t} \; c^4}{m^2 \; (m+2) \; \lambda^4} \left[ \frac{m}{2} \sum_{k=4}^{\infty} 
\frac{(\lambda t)^k}{k!} \; (k-3)(k+2) + \sum_{k=4}^{\infty} 
\frac{(\lambda t)^k}{k!} \; (k-3)(k-4) \right] \\ 
& \qquad (\text{see formulas (\ref{Lemma1Eq1}) and (\ref{Lemma2Eq1}) of Lemmas A1 and A2 
in Appendix below}) \\ 
& = \frac{8 \; e^{-\lambda t} \; c^4}{m^2 \; (m+2) \; \lambda^4} \biggl\{ \frac{m}{2} 
\bigl[ ((\lambda t)^2 -6) e^{\lambda t} + 2(\lambda t)^2 +6\lambda t + 6 \bigr] \\
& \hskip 4cm + ((\lambda t)^2 + 12) (e^{\lambda t} - 1) - 6\lambda t (e^{\lambda t} + 1) \biggr\} \\ 
& = \frac{8 \; c^4}{m^2 \; (m+2) \; \lambda^4} \biggl\{ \frac{m}{2} 
\bigl[ (\lambda t)^2 -6 + e^{-\lambda t} \bigl( 2(\lambda t)^2 +6\lambda t + 6 \bigr) \bigr] \\
& \hskip 4cm + ((\lambda t)^2 + 12) (1-e^{-\lambda t}) - 6\lambda t (1+e^{-\lambda t}) \biggr\} . 
\endaligned$$

Substituting this into (\ref{eq0191}) and taking into account that $(-i)^4 = 1$, after some simple calculations, we finally arrive at (\ref{SecondMomentFunct}). The theorem is thus completely proved. 
\end{proof}

\bigskip

One should emphasize that, from the symmetry of the process, it follows that 2-marginal second moment functions corresponding to any multi-index of the form $(0,\dots,0,2,0,\dots,0,2, 0,\dots,0)$) have the 
same form as in (\ref{SecondMomentFunct}). 

\bigskip

{\bf Remark 1.} Under the standard Kac scaling condition 
$$c\to\infty , \qquad \lambda\to\infty , \qquad \frac{c^2}{\lambda} \to \rho, \quad \rho>0,$$
function (\ref{SecondMomentFunct}) turns into: 

$$\aligned 
& \lim_{\substack{c, \; \lambda\to\infty\\(c^2/\lambda)\to\rho}} \; \mu_{(2,2,0,\dots,0)}(t) \\ 
& = \frac{8}{m^2 \; (m+2)} \; \lim_{\substack{c, \; \lambda\to\infty\\(c^2/\lambda)\to\rho}} 
\biggl\{ \left(\frac{c^2}{\lambda}\right)^2 
\; \left[ m \left( e^{-\lambda t} \left( t^2 + \frac{3t}{\lambda} + \frac{3}{\lambda^2} \right) 
+ \frac{t^2}{2} - \frac{3}{\lambda^2} \right) \right] \\ 
& \hskip 4cm + \biggl( t^2 +\frac{12}{\lambda^2} \biggr) \biggl( 1 - e^{-\lambda t} \biggr) 
+ \frac{6t}{\lambda} \biggl( 1 + e^{-\lambda t} \biggr) \biggr\} \\ 
& = \frac{8}{m^2 \; (m+2)} \; \rho^2 \; \left[ m \; \frac{t^2}{2} +t^2 \right] \\
& = \frac{4 \rho^2 t^2}{m^2} \\ 
& = \left( \frac{2 \rho t}{m} \right)^2,
\endaligned$$
and this is exactly the product of the variances of two coordinates of the $m$-dimensional homogeneous Brownian motion with zero drift and diffusion coefficient $\sigma^2 = 2\rho/m$ (we remind the fact 
that each coordinate of a multidimensional Brownian motion is an independent one-dimensional 
Brownian motion).  

\bigskip

{\bf Remark 2.} In the important three-dimensional case, the 2-marginal second moment function 
(\ref{SecondMomentFunct}) takes the form: 
\begin{equation}\label{SecondMomentFunct3} 
\mu_{(2,2,0)}(t) = \frac{8}{45} \; 
\frac{c^4}{\lambda^4} \biggl[ e^{-\lambda t} (2\lambda^2t^2 +3\lambda t - 3) + 
\frac{5}{2} \; \lambda^2t^2 - 6\lambda t +3 \biggr] . 
\end{equation}

The shape of moment function (\ref{SecondMomentFunct3}) is plotted in Fig. 1 below. 

\begin{center}
\begin{figure}[htbp]
\centerline{\includegraphics[width=10cm,height=8cm]{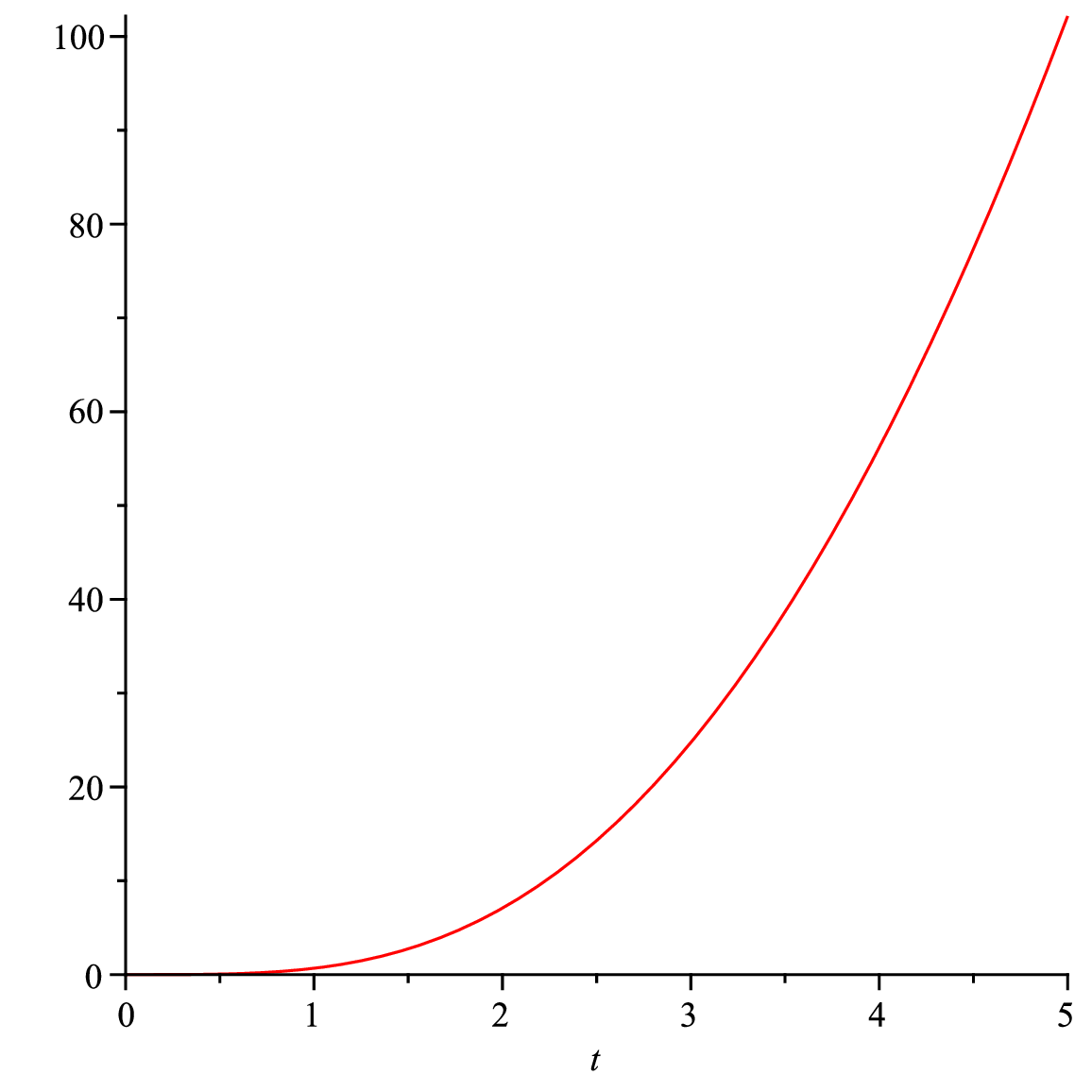}}
\caption{\it The shape of 2-marginal second moment function (\ref{SecondMomentFunct3}) on the interval 
$t\in [0, 5]$ (for $\lambda=1, \; c=2$)} 
\label{FigMomentFunction}
\end{figure}
\end{center} 

We see that this is a nonlinear monotonously increasing function whose growth becomes more and more 
steepen, as time increases. For small values of time variable, this function takes small values too. However, as time $t$ increases, the 2-marginal second moment function (\ref{SecondMomentFunct3}) grows like 
$C_1 t^2 - C_2 t$, where $C_1, C_2$ are some constants. By analyzing the 2-marginal second moment function 
(\ref{SecondMomentFunct}), one can assert that it has a similar behaviour in arbitrary dimension.

\bigskip 

\section{Conclusions and final remarks} 

The main achievement of the article is the explicit formula (\ref{SecondMomentFunct}) for the 2-marginal second moment function $\mu_{(2,2,0,\dots,0)}(t)$ of the symmetric Markov random flight in the Euclidean space $\Bbb R^m$ of arbitrary dimension $m\ge 3$. This result provides clear prospects for further generalizations in this area of ​​research. In particular, the problem of obtaining a similar exact formula for the 3-marginal second moment function $\mu_{(2,2,2,0,\dots,0)}(t)$ is of great interest. This importance is determined by the fact that such a formula would enable to obtain, as its particular case for $m=3$, the exact relation for the full second moment function $\mu_{(2,2,2)}(t)$ of the three-dimensional symmetric Markov random flight. This result, in totality with already known 1- and 2-marginal second moment functions 
(\ref{eq2910}) and (\ref{SecondMomentFunct3}), would give the complete solution to the problem of finding the second moment function and its marginals of the symmetric Markov random flight in the three-dimensional Euclidean space $\Bbb R^3$, which is of a special importance for stochastic analysis, statistical and quantum physics. 

To find a formula for the 3-marginal second moment function $\mu_{(2,2,2,0,\dots,0)}(t)$, one needs to use the method developed in this article and in \cite{kol0} and to try to derive an exact formula, similarly to  those presented by Propositions 1 and 2, for the coefficients $K_{\Vert\boldsymbol\alpha\Vert^6}(\gamma_n)$ 
at the term $\Vert\boldsymbol\alpha\Vert^6$ in the polynomials 
$\gamma_n=\gamma_n(\Vert\boldsymbol\alpha\Vert)$ for arbitrary $n\ge 7$. This, however, is a fairly difficult analytical problem, which will become the subject of forthcoming research. 

One can try to obtain relations for other $l$-marginals ($l\ge 4$) of second moment function, however such 
a result would not be so important for practical applications. The same also concerns other (the third, fourth, etc.) moment functions of the multidimensional Markov random flights.

\bigskip 

\begin{center}
{\bf Appendix} 
\end{center} 

In this appendix we prove two auxiliary lemmas that have been used in the proof of Theorem 1. 

\bigskip 

{\bf Lemma A1.} {\it For arbitrary real $x\in\Bbb R^1$, the following relation holds:}
\begin{equation}\label{Lemma1Eq1}
\sum_{k=4}^{\infty} \frac{x^k}{k!} \; (k-3) (k+2) = (x^2-6)e^x + 2 x^2 +6 x + 6 .
\end{equation} 
\begin{proof} 
We have: 

$$\aligned 
\sum_{k=4}^{\infty} \frac{x^k}{k!} \; (k-3) (k+2) & = \sum_{k=4}^{\infty} \frac{x^k}{k!} \; [k(k-1) - 6] \\
& = \sum_{k=4}^{\infty} \frac{x^k}{(k-2)!} - 6 \sum_{k=4}^{\infty} \frac{x^k}{k!} \\ 
& = \sum_{k=2}^{\infty} \frac{x^{k+2}}{k!} - 6 \sum_{k=4}^{\infty} \frac{x^k}{k!} \\ 
& = x^2 \; \sum_{k=2}^{\infty} \frac{x^k}{k!} - 6 \sum_{k=4}^{\infty} \frac{x^k}{k!} \\ 
& = x^2 \left( e^x -1 - x \right) - 6 \left( e^x -1 - x - \frac{x^2}{2!} - \frac{x^3}{3!} \right) \\ 
& = x^2 e^x - 6 e^x + 2 x^2 +6 x + 6 , 
\endaligned$$ 
proving (\ref{Lemma1Eq1}). 
\end{proof}

\bigskip 

{\bf Lemma A2.} {\it For arbitrary real $x\in\Bbb R^1$, the following relation holds:}
\begin{equation}\label{Lemma2Eq1}
\sum_{k=4}^{\infty} \frac{x^k}{k!} \; (k-3) (k-4) = (x^2 + 12) (e^x - 1) - 6x (e^x + 1) .
\end{equation} 

\begin{proof} 
We have: 
$$\aligned 
\sum_{k=4}^{\infty} \frac{x^k}{k!} \; (k-3) (k-4) & = \sum_{k=4}^{\infty} \frac{x^k}{k!} \; k(k-7) 
+ 12 \sum_{k=4}^{\infty} \frac{x^k}{k!} \\ 
& = \sum_{k=4}^{\infty} \frac{x^k}{(k-1)!} \; (k-7) + 12 \sum_{k=4}^{\infty} \frac{x^k}{k!} \\ 
& = x \; \sum_{k=3}^{\infty} \frac{x^k}{k!} \; (k-6) + 12 \sum_{k=4}^{\infty} \frac{x^k}{k!} \\ 
& = x \left[ \sum_{k=3}^{\infty} \frac{x^k}{(k-1)!} - 6 \sum_{k=3}^{\infty} \frac{x^k}{k!} \right] 
+ 12 \sum_{k=4}^{\infty} \frac{x^k}{k!} \\ 
& = x \left[ x \; \sum_{k=2}^{\infty} \frac{x^k}{k!} - 6 \sum_{k=3}^{\infty} \frac{x^k}{k!} \right] 
+ 12 \sum_{k=4}^{\infty} \frac{x^k}{k!} \\ 
& = x^2 \; \sum_{k=2}^{\infty} \frac{x^k}{k!} - 6 x \; \sum_{k=3}^{\infty} \frac{x^k}{k!}  
+ 12 \sum_{k=4}^{\infty} \frac{x^k}{k!} \\ 
& = x^2 (e^x -1 - x) - 6x \left( e^x -1 - x - \frac{x^2}{2!} \right) \\ 
& \qquad\qquad\quad + 12 \left( e^x -1 - x - \frac{x^2}{2!} - \frac{x^3}{3!} \right) \\ 
& = (x^2 + 12) (e^x - 1) - 6x (e^x + 1) , 
\endaligned$$ 
proving (\ref{Lemma2Eq1}). 
\end{proof}

\bigskip 

{\bf Declaration.} The author declares no potential conflicts of interest with respect to the research, authorship, and/or publication of this article. The author has no data availability to share.

\end{document}